\documentclass[10pt]{article}
\usepackage{amsmath, amsfonts, amsthm, amssymb,verbatim}
        \newtheorem{theorem}{Theorem}[section]

        \newtheorem{lemma}[theorem]{Lemma}

\numberwithin{equation}{section}

\newcommand{\N}{\mathbb{N}}
\newcommand{\R}{\mathbb{R}}
\newcommand{\Leb}{{\mathcal L}}
\newcommand{\Haus}{{\mathcal H}}
\newcommand{\res}{\mathop{\hbox{\vrule height 7pt width .5pt depth 0pt
\vrule height .5pt width 6pt depth 0pt}}\nolimits}
\newcommand{\loc}{{\rm loc}}

\newcommand{\spt}{{\rm supp \, }}
\newcommand{\lip}{{\rm Lip}}
\newcommand{\D}{\mathcal{D}}
\newcommand{\I}{\mathcal{I}}
\newcommand{\II}{\bar{\mathcal {I}}}

\renewcommand{\flat}{\mathbb{F}}
\newcommand{\mass}{\mathbb{M}}
\renewcommand{\d}{\mathbf{d}}
\renewcommand{\i}{\mathbf{i}}
\newcommand{\F}{\mathcal F}

\begin{document}

\title
{Leaf superposition property for integer rectifiable currents}
 
\author 
{
Luigi Ambrosio$^1$, Gianluca Crippa$^2$, 
\\
and  
\\  
Philippe G. LeFloch$^3$}
\date{January 15, 2008}
\maketitle

\footnotetext[1]{Scuola Normale Superiore, Piazza dei Cavalieri 7, 56126 Pisa, Italy. 
Email: l.ambrosio@sns.it.}

\footnotetext[2]{Dipartimento di Matematica, 
Universit\`a degli Studi di Parma, Viale G.P.~Usberti 53/A (Campus), 43100 Parma, Italy. 
E-mail: gianluca.crippa@unipr.it.}

\footnotetext[3]{Laboratoire Jacques-Louis Lions \& Centre National de la Recherche Scientifique,
Universit\'e Pierre et Marie Curie (Paris 6), 4 Place Jussieu,  75252 Paris, France.
E-mail : {pgLeFloch@gmail.com.} 
{\tt Cite this paper as:}  L. Ambrosio, G. Crippa, and P.G. LeFloch, 
Leaf superposition property for integer rectifiable currents, 
Netw. Heterog. Media 3 (2008), 85--95.}
 
\begin{abstract}
We consider the class of integer rectifiable currents without boundary in $\R^n\times\R$
satisfying a positivity condition.
We establish that these currents can be written as a linear superposition of graphs of
finitely many functions with bounded variation.
\end{abstract}

\section{Introduction and statement of the main result}

It is well known that a locally integrable function in $\R^n$
belongs to $BV_\loc$ (the space of functions of locally bounded variation)
if and only if its subgraph has locally finite perimeter in $\R^n\times\R$.
The connections between the analytic properties of $u$ and the geometric properties
of its (sub)graph are well described, using the more powerful language of currents,
in \cite[4.5.9]{Federer} or \cite[4.1.5]{GMS}.
Recall that currents provide a very natural setting to discuss analytic problems
with a geometrical content, and have been successfully used in
many areas. In particular, Giaquinta, Modica, and Sou\v{c}ek  introduced the
notion of Cartesian current and used it to attack many problems in
the calculus of variations (see the extensive monograph \cite{GMS})
including non-linear elasticity, harmonic maps between manifolds,
relaxed energies, etc.

The aim of this paper is to show the representation of a suitable class
of integer rectifiable currents in $\R^n\times\R$
as the superposition of finitely many graphs (referred to as ``leaves'')
of functions with bounded variation. In some sense this result has some
connections with Almgren's theory \cite{almgren1}, \cite{almgren2} (developed in arbitrary dimension and
codimension) of approximation, up to sets of small measure, of (minimal) currents
by multi-valued Lipschitz graphs: here the regularity condition is weakened to $BV$,
and this allows a complete description of the current, at least in codimension one, as a multi-valued graph.
We rely on techniques of geometric measure theory, especially the concept of $BV$ maps and currents in metric
spaces developed in Ambrosio~\cite{Ambrosio} and Ambrosio and Kirchheim~\cite{AK}.

We refer to the following section for the notation and state now the main result of this paper.
If $ u : \R^n \to \R$ is a locally $BV$ function, we
denote by $\i (u)$ the $n$-dimensional boundary-free current canonically associated with
the graph of $u$ in $\R^n \times \R$, obtained
(roughly speaking) by completion of the discontinuities of $u$ with vertical segments.

\begin{theorem}\label{t:main}
Let $T \in \I_n(\R^{n+1})$ be an $n$-dimensional integer rectifiable current in $\R^{n+1}=\R^n_x \times \R_y$
satisfying the zero-boundary condition $\partial T = 0$,
the positivity condition $T \res dx \geq 0$ and the cylindrical mass condition
\begin{equation}
\label{e:massass}
\mass_{B_R(0) \times \R} (T) < \infty \qquad \text{ for every $R>0$.}
\end{equation}
Then, there exist a unique integer $N$ and a unique family of functions\break $u_j \in BV_\loc(\R^n;\R)$,
${1 \leq j \leq N}$, satisfying
\begin{equation}\label{e:monot}
u_1 \leq u_2 \leq \ldots \leq u_N,
\end{equation}
such that the given current $T$ is the superposition of the canonical Cartesian currents $\i(u_j)$
associated with the functions $u_j$, that is,
\begin{equation}\label{e:super}
T = \sum_{j=1}^N \i (u_j).
\end{equation}
In addition, the following additivity property holds:
\begin{equation}\label{nocance}
\| T \| = \sum_{j=1}^N \| \i (u_j) \|.
\end{equation}
\end{theorem}

We call each function $u_j$ a leaf of the decomposition of $T$, and
we refer to \eqref{e:super} as the {\em canonical leaf
decomposition} of $T$. Heuristically \eqref{nocance} follows from
\eqref{e:super} because all graphs have a common orientation in
their intersection, so that no cancellations occur; notice that the
additivity property \eqref{nocance} does not hold for more general
decompositions which satisfy condition \eqref{e:super}, but not the
monotonicity assumption \eqref{e:monot}.

For an application of this result we refer to \cite{ACL}, where a
geometric approach to tackle multi-dimensional scalar conservation
laws is developed. Therein, solutions are defined geometrically as
currents, rather than as functions satisfying entropy inequalities.
The leaf decomposition is used to show the existence of entropy solutions in
this setting, as the superposition of graphs of entropy
solutions. (See \cite{ACL} for details.)

\section{Preliminaries and notation}

\subsection{Currents}

We denote by $\D_m(\R^k)$ the space of $m$-dimensional currents in
$\R^k$, that is the dual space of all linear and continuous
functionals defined on the space $\D^m(\R^k)$ of all smooth and
compactly supported differential $m$-forms. The space $\D_m(\R^k)$
is equipped with the usual weak-star topology induced by this
duality. The duality bracket between a current $T \in \D_m(\R^k)$
and a form $\omega \in \D^m(\R^k)$ is denoted by $\langle T, \omega
\rangle$.

The boundary of a current $T \in \D_m(\R^k)$ is the current $\partial T \in \D_{m-1}(\R^k)$ defined by
$$
\langle \partial T, \omega \rangle
= \langle T, d \omega \rangle, \qquad \omega \in \D^{m-1}(\R^k),
$$
where $d \omega \in \D^m(\R^k)$ denotes the differential of the form
$\omega \in \D^{m-1}(\R^k)$. If $T \in \D_m(\R^k)$ and $\alpha \in
\D^h(\R^k)$ for some $h \leq m$, we denote by $T \res \alpha \in
\D_{m-h}(\R^k)$ the saturation of the current $T$ with the form
$\alpha$, which is defined by
$$
\langle T \res \alpha, \omega \rangle
= \langle T, \alpha \wedge \omega \rangle,
\qquad  \omega \in \D^{m-h}(\R^k).
$$

The (local) mass of a current $T \in \D_m(\R^k)$ is defined for every open set $\Omega \subset \R^k$ as
$$
\mass_\Omega (T) = \sup \left\{ \langle T, \omega \rangle \; : \;
\omega \in \D^m(\R^k),\; \spt \omega \subset \Omega,\; \|\omega\| \leq 1 \right\}.
$$
If $T$ has locally finite mass, the set function $\Omega\mapsto\mass_\Omega(T)$ is the
restriction to bounded open sets of a nonnegative Radon measure that we shall denote by
$\|T\|$, so that $\|T\|(\Omega)=\mass_\Omega(T)$ for all bounded open sets
$\Omega\subset\R^k$.
Given a current $T \in \D_m(\R^k)$ with locally finite mass, there exists a unique
(up to $\|T\|$-negligible sets) $\| T \|$-measurable map $\vec{T}$ defined
on $\R^k$ and with values in the set of $m$-vectors such that $\vec{T}$ is a unit
$m$-vector $\| T \|$-almost everywhere and (here $\langle \cdot,\cdot\rangle$ is the
standard duality between $m$-vectors and $m$-covectors)
\begin{equation}\label{may1}
\langle T, \omega \rangle = \int_{\R^k} \langle \vec{T}(x),
\omega(x) \rangle \, d \|T\|(x),
 \qquad
 \omega \in \D^m(\R^k).
\end{equation}
Whenever \eqref{may1} holds, we shall write $T=\vec{T}\|T\|$.

We will be especially interested in the subclass $\I_m(\R^k) \subset
\D_m(\R^k)$ of all $m$-dimensional integer rectifiable currents $T$
for which, by definition, there exists a triple $(M, \theta,\tau)$,
where $M\subset\R^k$ is a countably $\Haus^m$-rectifiable set,
$\theta : M \to \N\setminus\{0\}$ is a locally integrable function
and $\tau$ is a Borel orientation of $M$ (i.e. a Borel map
$x\mapsto\tau(x)=\xi_1(x)\wedge\ldots\wedge\xi_m(x)$ with values in
unit and simple $m$-vectors whose span is the approximate tangent
space to $M$ at $x$) such that $T=\tau\theta\Haus^m\res M$, or
equivalently $\vec{T}=\tau$ and $\|T\|=\theta\Haus^m\res M$. We
shall also write $T=(M,\theta,\tau)$, and we refer to $M$ as the
support of $T$ and to $\theta$ as the multiplicity of $T$ (both are
uniquely determined up to $\Haus^m$-negligible sets).

\subsection{$0$-dimensional integer rectifiable currents with finite mass}

In this section we consider a very special class of integer rectifiable
currents, the $0$-dimensional ones with finite mass on the real line
$\R$. We denote by $\II_0(\R)$ the set of these currents and we
notice that it consists
of those currents that can be expressed as a finite sum of Dirac masses with
weight $\pm 1$. This means that every $S \in \II_0(\R)$ can be
written as
$$ S = \sum_{j=1}^l \sigma_j \delta_{A_j},
$$
where the $A_j$ are (not necessarily distinct) points of $\R$ and
$\sigma_j = \pm 1$. We will call {\em average} of the current $S \in
\II_0(\R)$ the integer $\sum_j \sigma_j$. For every $h \in \N$ we
denote by $\I_0^h (\R)\subset\II_0(\R)$ the set consisting of all
nonnegative $0$-dimensional integer rectifiable currents in $\R$
with average $h$:
\begin{equation}\label{e:uff}
\I_0^h (\R) := \Big\{ S \in \II_0(\R) \; : \; S = \sum_{j=1}^h
\delta_{A_j} \Big\} ;
\end{equation}
notice again that the points $A_j \in \R$ need not be distinct.

On the set $\II_0(\R)$ we define
$$
\flat (S) := \sup \left\{ \langle S, \phi \rangle \; : \; \phi \in
\lip_{b,1}(\R) \right\},
 \qquad S \in \II_0(\R),
$$
where $\lip_{b,1}(\R)$ denotes the set of bounded real-valued
Lipschitz functions defined on $\R$ with Lipschitz constant less or
equal than one. Notice that, if $S \in \II_0(\R)$ has non-zero
average, then obviously $\flat (S) = + \infty$; on the other hand
$$\flat (S)\leq \mass_\R(S) \; {\rm diam\,}({\rm supp\,}S)<+\infty$$
for all $S\in
\II_0(\R)$ with zero average. It is also immediate to check that, for $S = \delta_A -
\delta_B$, we have $\flat (S) = |A-B|$. A generalization of this
fact is given by the following well-known lemma.

\begin{lemma}\label{l:kant}
If $S$ and $S' \in \II_0(\R)$ are of the form
$$
S = \sum_{j=1}^h \delta_{A_j}, \qquad S' = \sum_{j=1}^h \delta_{B_j},
$$
with $A_1 \leq A_2 \leq \ldots \leq A_h$ and $B_1 \leq B_2 \leq
\ldots \leq B_h$, then
\begin{equation}\label{gattini}
\sum_{j=1}^h | A_j - B_j| = \flat ( S-S').
\end{equation}
\end{lemma}

\begin{proof} We give an elementary proof, which uses ideas from the
theory of optimal transportation (see \cite{villani}). We notice first that
the inequality $\geq$ in \eqref{gattini} is an obvious consequence
of the inequality $|A_j-B_j|\geq|\phi(A_j)-\phi(B_j)|$ for all
$\phi\in \lip_{1,b}(\R)$, so we need only to build $\phi\in
\lip_{1,b}(\R)$ such that
\begin{equation}\label{gattini1}
\sum_{j=1}^h | A_j - B_j| \leq \langle S-S',\phi\rangle.
\end{equation}
By the compactness of the support of $S-S'$, it suffices to
construct a 1-Lipschitz function $\phi$ with this property. To this
aim, we first notice that the fact that the list of the $A_j$'s and
of the $B_j$'s are ordered implies
\begin{equation}\label{gattini2}
\sum_{j=1}^h |A_j-B_j|\leq\sum_{j=1}^h |A_j-B_{\sigma(j)}|
\end{equation}
for any permutation $\sigma$ of $\{1,\ldots,h\}$ (this can be seen
by showing that the right hand side does not increase if a
permutation $\sigma$ with $B_{\sigma(i)}> B_{\sigma(j)}$ for some
$i<j$ is replaced by another one $\tilde\sigma$ with
$\tilde\sigma(i)=\sigma(j)$, $\tilde\sigma(j)=\sigma(i)$ and
$\tilde\sigma(k)=\sigma(k)$ for $k\neq i,\,j$). More generally, one
can use \eqref{gattini2} and the fact that permutation matrices are
extremal points in the class of bi-stochastic matrices to obtain
(the so-called Birkhoff theorem, see \cite{villani})
\begin{equation}\label{gattini6}
\sum_{j=1}^h |A_j-B_j|\leq\sum_{i,j=1}^h m_{ij}|A_j-B_i|
\end{equation}
for any nonnegative $m_{ij}$ with $\sum_i m_{ij}=\sum_i m_{ji}=1$
for all $j=1,\ldots,h$.

The minimization of the functional
$m\mapsto\sum_{i,j}m_{ij}|A_j-B_i|$ subject to the above constraints
on $m$ is a (very) particular case of Monge-Kantorovich optimal
transport problem of finding an optimal coupling between $S$ and
$S'$ with cost function $c(x,y)=|x-y|$. Kantorovich's duality theory
gives that the infimum of this problem, namely $\sum_j|A_j-B_j|$, is
(see \cite{villani} again, where an explicit construction of the
maximizing $\phi$ is given)
$$
\max_{\phi\in\lip_1(\R)}\langle S-S',\phi\rangle.
$$
\end{proof}

For every fixed $h\in \N$ we define
$$
\d (S,S') := \flat (S-S') = \sup\left\{ \langle S, \phi \rangle -
\langle S', \phi \rangle \; : \; \phi \in \lip_{b,1}(\R) \right\},
 \quad  S, \,S' \in \I_0^h(\R),
$$
which is easily seen to be a finite distance in $\I_0^h(\R)$ (indeed,
since $S$ and $S'$ belong to the same set $\I_0^h(\R)$, the difference $S - S'$ has zero average).

\subsection{Slices of a current}

Given $T \in \I_n(\R^n\times\R)$ we consider the vertical slices of $T$ at $x\in\R^n$,
$$
T_x := \langle T, dx, x \rangle \in \I_0(\R),
$$
see for instance \cite{AK}, \cite{simon}. This family of currents is
uniquely determined, up to $\Leb^{n}$-negligible sets, by the
identity $\int_{\R^n} T_x\,dx=T\res dx$, i.e.
\begin{equation}\label{l_slice}
\int_{\R^n} \langle T_x, \varphi(x,\cdot) \rangle \, dx =
\langle T\res dx,\varphi\rangle
\end{equation}
for all $\varphi\in C^\infty_c(\R^n\times\R)$.
Furthermore, the masses of $T_x$ are related to the mass of $T$ by
\begin{equation}\label{cylindr}
\int_\Omega \mass_{\R}(T_x)\,dx\leq \mass_{\Omega\times\R}(T)
\end{equation}
for all bounded open sets $\Omega\subset\R^n$. As a consequence,
$T_x\in\II_0(\R)$ for $\Leb^n$-a.e.~$x\in\Omega$ whenever
$\mass_{\Omega\times\R}(T)<+\infty$.

\subsection{The current associated to the graph of a $BV$ function}

Recall that $u\in L^1_\loc(\R^n)$ is said to be a locally $BV$ function if its
distributional derivative $Du=(D_1u,\ldots,D_nu)$ is an $\R^n$-valued
measure with locally finite total variation in $\R^n$, and we shall denote by
$\|Du\|$ this total variation.

In this section we are going to describe how we can canonically associate to
$u\in BV_\loc(\R^n)$ a current $\i(u)\in\I_n(\R^n\times\R)$ with no boundary,
finite mass on cylinders $\Omega\times\R$ with $\Omega$ bounded, and
satisfying
\begin{equation}\label{carto}
\langle \i(u),\varphi dx\rangle:=\int_\Omega \varphi(x,u(x))\,dx
\qquad\forall\varphi\in C^\infty_c(\R^n\times\R).
\end{equation}
These two conditions are actually sufficient to characterize a unique current,
see step 5 of the proof of Theorem~\ref{t:main}.

Geometrically, this current corresponds to the integration on the graph of
$u$, with the orientation induced by the map $x\mapsto (x,u(x))$, and this
description works perfectly well when $u\in C^1$. In order to
define $\i(u)$ in the general case when $u\in BV_\loc$, we first
define the subgraph $E(u)$ of $u$ by
$$
E(u):=\left\{(x,y)\in \R^n\times\R:\ y\leq u(x)\right\}.
$$
It is well known that $E(u)$ has locally finite perimeter in
$\R^n\times\R$ (i.e. $\chi_{E(u)}\in BV_\loc(\R^n\times\R)$), so it
has a measure-theoretic boundary (the set of points where the
density of $E(u)$ is neither 0 nor 1), that we shall denote by
$\Gamma(u)$. De Giorgi's theorem on sets of finite perimeter ensures
that $\Gamma(u)$ is countably $\Haus^n$-rectifiable, and that
\begin{equation}\label{basic}
D\chi_{E(u)}=-\nu_{E(u)}\Haus^n\res\Gamma(u)
\end{equation}
(the unit vector $\nu_{E(u)}$ is the so-called approximate outer
normal to $E(u)$). Then, we define
\begin{equation}\label{gat5}
\i(u):=(\Gamma(u),1,\tau_u),
\end{equation}
where $\tau_u$ is the unit $n$-vector spanning $\nu_{E(u)}^\perp$ (the approximate
tangent space to $\Gamma(u)$), characterized by
$$
\langle  dx_1\wedge\cdots\wedge dx_n\wedge dy \; , \; \tau_u\wedge\nu_{E(u)} \rangle \geq 0 .
$$
Equivalently, invoking the relation \eqref{basic}, we can define
$$
\aligned \langle \i(u),\varphi dx\rangle & := - \int_{\R^n\times\R}\varphi
\,dD_y \chi_{E(u)}, \\
\langle \i(u),\varphi\widehat{dx}_j\wedge
dy\rangle & :=\int_{\R^n\times\R}\varphi\,dD_j\chi_{E(u)}, \quad
j=1,\ldots,n
\endaligned
$$
(here $\widehat{dx}_j:=(-1)^{n-j}dx^1\wedge\cdots\wedge
dx^{j-1}\wedge dx^{j+1}\wedge\cdots\wedge dx^n$). In the case $u\in
C^1(\R^n)$, using the area formula, it is easy to check that this
definition coincides with the geometric picture, and in particular
that $\partial \big( \i(u) \big) = 0$ and \eqref{carto} hold. In the general case both can
be obtained, for instance, by approximation (notice that $u_i\to u$
in $L^1_\loc$ implies $E(u_i)\to E(u)$ in $L^1_\loc$ and therefore
weak convergence of the associated currents).

We will need the following strong locality property of $\tau_u$.

\begin{lemma}\label{localtau}
Let $u,\,v\in BV_\loc(\R^n)$ with $u\geq v$. Then $\tau_u=\tau_v$ $\Haus^n$-a.e.~on
$\Gamma(u)\cap\Gamma(v)$.
\end{lemma}
\begin{proof}
It suffices to show that $\nu_{E(u)}=\nu_{E(v)}$ $\Haus^n$-a.e.~on $\Gamma(u)\cap\Gamma(v)$.
It is a general property of sets of finite perimeter $E\subset\R^{n+1}$ that, for $\Haus^n$-a.e.
$w\in\partial^* E$, the rescaled sets $(E-w)/r$ converge in $L^1_\loc$ as $r\downarrow 0$
to the halfspace having $\nu_E(w)$ as outer normal. In our case,
$E(u)\supset E(v)$ because $u\geq v$, so that all points $w$ where both $(E(u)-w)/r$
and $(E(v)-w)/r$ converge to a halfspace, the halfspace has to be the same. This implies
the stated equality $\Haus^n$-a.e.~of the outer normals.
\end{proof}

\subsection{Metric spaces valued $BV$ functions}\label{ss:mbv}

We now recall the main features of the theory of $BV$ functions with
values in a metric space, developed in Ambrosio \cite{Ambrosio} and
Ambrosio and Kirchheim \cite{AK}. Let $(E,d)$ be a metric space such
that there exists a countable family $\F \subset \lip_{b,1}(E)$
which generates the distance, in the sense that
$$
d(x,y) = \sup_{\Phi \in \F} | \Phi(x) - \Phi(y)|, \qquad  x,\,y \in E.
$$
We say that a function $f : \R^n \to E$ is a function of metric locally bounded variation,
and we write $f \in MBV_\loc(\R^n;E)$, if $\Phi \circ f \in BV_\loc(\R^n)$ for every $\Phi \in \F$
and if there exists a positive locally finite measure $\nu$ in $\R^n$ such that
$$
\nu \geq \left\| D ( \Phi \circ f) \right\|, \qquad  \Phi \in \F.
$$
The minimal $\nu$ such that the previous condition holds will still
be denoted by $\| Df\|$. It is possible to check that the class
$MBV_\loc(\R^n;E)$ and the measure $\| Df \|$ are independent of the
choice of the family $\F$.

We now consider the metric space $(\I_0^h(\R),\d)$ previously
defined. To every $\phi \in \lip_{b,1} (\R)$ we associate the map
$\Phi_\phi \in \lip_1(\I_0^h(\R))$ defined by
$$
\Phi_\phi (S) := \langle S, \phi\rangle, \qquad S \in \I_0^h(\R).
$$
Indeed, it is immediate to check the Lipschitz continuity
$$
| \Phi_\phi(S) - \Phi_\phi(S') | = | \langle S, \phi \rangle -
\langle S', \phi \rangle | \leq \d (S, S'), \qquad S,\,S' \in
\I_0^h(\R).
$$
By a standard density argument, it is possible to select a countable
family $\F\subset \lip_{b,1}(\R)\cap C^\infty(\R)$ with the property that
\begin{equation}\label{gat3}
\d (S,S') = \sup_{\phi\in\F}
\left\{ \langle S, \phi\rangle - \langle S', \phi\rangle \right\},
 \qquad S,\,S' \in \I_0^h(\R).
\end{equation}

\begin{lemma}\label{easy}
Let $E$ and $F$ be metric spaces. Then $M\circ f\in MBV_\loc(\R^k;F)$ whenever $f\in MBV_\loc(\R^k;E)$ and
$M:E\to F$ is an $L$-Lipschitz function, and $\|D(M\circ f)\|\leq
L\|D f\|$. Furthermore, $MBV_\loc(\R^k;\R)$ coincides with
$BV_\loc(\R^k)$.
\end{lemma}
\begin{proof}
Let $\phi\in\lip_{b,1}(F)$, $g=M\circ f$ and $\psi=\phi\circ M$; then $\psi\in \lip_b(E)$
and its Lipschitz constant is less than $L$; as a consequence, $\|D ( \psi\circ f ) \|\leq L\|Df\|$.
Since $\psi\circ f=\phi\circ g$ we obtain that $g\in MBV_\loc(\R^k;F)$ and $\|D g\|\leq L\|Df\|$.\\
The inclusion $BV_\loc(\R^k)\subset MBV_\loc(\R^k;\R)$ is a simple consequence of the
stability of $BV$ functions under left composition with Lipschitz maps; to prove the
opposite inclusion, let $f\in MBV_\loc(\R^k;\R)$ and fix an open ball $B\subset\R^k$;
by definition all truncated functions $f_a:=-a\lor (f\land a)$ belong to $BV(B)$ and
$\|Df_a\|\leq \|Df\|$, since we can see $f_a$ as the composition of $f$ with the
map $\eta_a \in \lip_{b,1}(\R)$ defined as the identity for $x \in [-a,a]$, as the constant $a$
for $x > a$ and as the constant $-a$ for $x < -a$. Therefore, denoting by $\bar f_a$
their averages in $B$,
by Poincar\'e inequality we obtain that $f_a-\bar f_a$ is bounded in $L^1(B)$.
Thanks to the compactness of the embedding of $BV$ in $L^1$, we can find a sequence
 $a_i\to +\infty$ such that $\bar f_{a_i}$ converges to some
$m\in\overline\R$ and $f_{a_i}-\bar f_{a_i}$ converge in $L^1(B)$
and $\Leb^n$-almost everywhere to $g\in BV(B)$: if $m\in\R$ we
immediately obtain that $f=m+g\in BV(B)$. If not, we obtain that
$|f|=+\infty$ $\Leb^n$-almost everywhere, contradicting the
assumption that $f$ is real valued.
\end{proof}

\section{Proof of the main theorem}

This section is entirely devoted to the proof of Theorem
\ref{t:main}. We address separately the existence of the
decomposition, its uniqueness and the equality of the total
variations. In the course of the proof we will occasionally use forms $\omega$
in $\R^n\times\R$ whose supports are not compact, but have a compact
projection on $\R^n$. Their use can be easily justified by a
truncation argument, based on the fact that the currents under
consideration have finite mass on cylinders $\Omega\times\R$ with
$\Omega\subset\R^n$ bounded.

\subsection{Existence of a decomposition} We proceed in 5 steps.

{\bf Step 1.} We begin by proving that there exists an integer $N$
(depending on $T$ only) such that, for $\Leb^n$-a.e.~$x \in
\R^n$, the slice $T_x \in \II_0(\R)$ is the sum of $N$ Dirac masses
with unit weight: more precisely, for $\Leb^n$-a.e.~$x \in
\R^n$ there exist $N$ real values
\begin{equation}\label{e:ord}
u_1(x) \leq u_2(x) \leq \ldots \leq u_N(x)
\end{equation}
satisfying
\begin{equation}\label{e:struct}
T_x = \sum_{j=1}^N \delta_{u_j(x)}.
\end{equation}

We first show that $T_x \geq 0$. Fix two nonnegative functions $\phi \in C^\infty_c(\R)$
and $\psi \in C^\infty_c(\R^n)$, and apply \eqref{l_slice} with $\varphi(x,y)=\psi(x)\phi(y)$ to get
$$
\int_{\R^n} \langle T_x, \phi \rangle\, \psi(x) \, dx
= \langle T \res dx,\varphi \rangle \geq 0,
$$
since we assumed $T \res dx \geq 0$. Hence, by the arbitrariness of
$\psi$, we deduce that $\langle T_x, \phi\rangle \geq 0$ for
$\Leb^n$-a.e.~$x \in \R^n$. By a simple density argument we can
obtain an $\Leb^n$-negligible set $E$ independent of $\phi$ such that
$\langle T_x,\phi\rangle \geq 0$ for all $\phi\in C^\infty_c(\R)$
and $x\in\R^n\setminus E$. This proves that $T_x\geq 0$ for all
$x\in\R^n\setminus E$.

Knowing that $T_x \geq 0$, the mass of $T_x$ is simply given by
$\langle T_x, 1 \rangle$ (notice that this function is locally
integrable by \eqref{cylindr} and assumption \eqref{e:massass}, and
takes $\Leb^{n}$-almost everywhere its values in $\N$ because
$\Leb^{n}$-almost all the slices are integer rectifiable).

We want to show that the map $x \mapsto \langle T_x, 1 \rangle$ is $\Leb^n$-equivalent to a
constant in $\R^n$.  Indeed, for every function $\psi \in C^\infty_c(\R^n)$ we can compute
(applying again \eqref{l_slice})
$$
\int_{\R^n} \langle T_x, 1 \rangle \frac{\partial \psi}{\partial x_i}(x) \, dx
= \langle T \res dx,\frac{\partial \psi}{\partial x_i}\rangle
= (-1)^{n-1} \langle T, d \left( \psi \widehat{dx}_i \right) \rangle = 0,
$$
since $\partial T = 0$. Hence we denote by $N\in\N$ the
$\Leb^n$-a.e.~constant value of $\langle T_x, 1 \rangle$, and we can
obviously assume that $N\geq 1$. In view of the representation
\eqref{e:uff}, this means that $T_x \in \I_0^N(\R)$ for
$\Leb^n$-a.e.~$x \in \R^n$. This leads us to the decomposition
\eqref{e:ord}--\eqref{e:struct}.

\

{\bf Step 2.} Next, we claim that the map
$$
\aligned & \R^n \to (\I_0^N(\R), \d),
\\
& x \mapsto T_x,
\endaligned
$$
belongs to $MBV_\loc(\R^n;\I_0^N(\R))$.

We proceed as in the proof of Theorem 8.1 of \cite{AK}. Recalling
the definitions and the discussion in Subsection \ref{ss:mbv}, we
only need to show that for every $\phi \in \lip_{b,1}(\R)\cap
C^\infty(\R)$ the map
$$ x \mapsto \langle T_x, \phi \rangle $$
belongs to $BV_\loc(\R^n)$, with a uniform (with respect to $\phi$) control of the derivative.

For every $\psi \in C^\infty_c(\R^n)$, applying once more \eqref{l_slice} we compute
$$
\int_{\R^n} \langle T_x, \phi \rangle \frac{\partial \psi}{\partial x_i} (x) \, dx
= \langle T \res dx ,\frac{\partial \psi}{\partial x_i}  \phi \rangle
= \langle T, \frac{\partial \psi}{\partial x_i}\phi \, dx \rangle
= - \langle T, \phi'\psi \, \widehat{dx}_i\wedge dy \rangle,
$$
using in the last equality the fact that $ \partial T = 0$.
Therefore, taking the modulus of both sides, we obtain
$$
\left|\int_{\R^n} \langle T_x, \phi \rangle \frac{\partial \psi}{\partial x_i} (x) \, dx \right|
\leq \int_{\R^n}|\psi|\,d\pi_\#\|T\|,
$$
where $\pi:\R^n\times\R\to\R^n$ is the projection on the $x$ variable.
This implies that the total variation of the distributional derivative of
$x\mapsto\langle T_x,\phi\rangle$ satisfies
$$
\left\| D\langle T_x, \phi \rangle \right\| \leq n \pi_\#\| T \|.
$$

\

{\bf Step 3.} Given $S \in \I_0^N(\R)$ of the form
$$
S = \sum_{j=1}^N \delta_{A_j}, \qquad \text{ with } A_1 \leq A_2 \leq \ldots \leq A_N,
$$
let us prove that the map
$$
\aligned
& (\I_0^N(\R),\d) \to \R,  \\
& S \mapsto A_N,
\endaligned
$$
is $1$-Lipschitz continuous.

Let $S \in \I_0^N(\R)$ be of the form above and $S' \in
\I_0^N(\R)$ be of the same form
$$
S' = \sum_{j=1}^N \delta_{A_j'}, \qquad \text{ with } A_1' \leq A_2' \leq \ldots \leq A_N'.
$$
Then
$$
| A_N - A_N' | \leq \sum_{j=1}^N | A_j - A_j' | = \flat \left( \sum_{j=1}^N \delta_{A_j} - \sum_{j=1}^N \delta_{A_j'} \right)
= \d (S,S'),
$$
where we have used Lemma \ref{l:kant}.

\

{\bf Step 4.} Finally we claim that the map
$$ x \mapsto u_N(x), \qquad \R^n \to \R$$
belongs to $BV_\loc(\R^n)$.

We have already seen in Step 2 that the map
$$ x \mapsto T_x, \qquad \R^n \to \I_0^N(\R)$$
is $MBV_\loc$ and in Step 3 that the map defined by
$$ \sum_{j=1}^N \delta_{z_j} \mapsto \max_{1\leq i\leq n}z_i, \qquad
\I_0^N(\R) \to \R $$ is Lipschitz continuous. Then, Lemma~\ref{easy} yields
that their composition, namely $u_N$, belongs to
$MBV_\loc(\R^n;\R)$, which is nothing but $BV_\loc(\R^n)$.

\

{\bf Step 5.} Induction and conclusion of the proof.

Up to now we have selected the top leaf of the decomposition. Now define
$$ \hat{T} = T - \i(u_N) . $$
It is readily checked that $\hat{T}$ is an $n$-dimensional integer
rectifiable current in $\R^{n+1}$, satisfying the zero-boundary
condition, the positivity condition and the cylindrical mass
condition as in the statement of the theorem, and that for
$\Leb^n$-a.e.~$x \in \R^n$ we have
$$ \hat{T}_x = T_x - \delta_{u_N(x)} = \sum_{j=1}^{N-1} \delta_{u_j(x)} . $$
Then, it suffices to apply again $N-1$ times the construction
described in the previous steps to deduce that all functions $u_j$
belong to $BV_\loc(\R^n)$ and, by construction,
$\left(T-\sum_{j=1}^N \i(u_j)\right)\res dx=0$. Let now
$R:=T-\sum_{j=1}^N \i(u_j)$ and let us prove that $\partial R=0$ and
$R\res dx=0$ imply $R=0$. Indeed, given $\psi\in
C^\infty_c(\R^n\times\R)$, let
$\varphi(x,y):=\int_{-\infty}^y\psi(x,s)\,ds$; then for every $j=1,\ldots,n$ we have
$$
0=\langle\partial R,\varphi \widehat{dx}_j\rangle=
(-1)^{n-1} \langle R,\frac{\partial\varphi}{\partial x_j}dx+\psi \widehat{dx}_j \wedge dy \rangle=
(-1)^{n-1} \langle R,\psi  \widehat{dx}_j \wedge dy \rangle.
$$

Finally, property \eqref{e:monot} is a consequence of the choice we
have done in \eqref{e:ord}.

\subsection{Uniqueness of the decomposition} The uniqueness of this decomposition is immediate.
Assume that we have two decompositions
$$ T = \sum_{j=1}^N \i (u_j) = \sum_{j=1}^M \i (v_j) ,$$
with $u_j \in BV_\loc(\R^n)$ for $j=1,\ldots,N$ and $v_j \in BV_\loc(\R^n)$ for $j=1,\ldots,M$ satisfying
\begin{equation}\label{e:dord}
u_1 \leq u_2 \leq \ldots \leq u_N \quad \text{ and } \quad v_1 \leq v_2 \leq \ldots \leq v_M .
\end{equation}
For $\Leb^n$-a.e.~$x \in \R^n$ the slice $T_x$ satisfies
$$ T_x = \sum_{j=1}^N \delta_{u_j(x)} = \sum_{j=1}^M \delta_{v_j(x)} . $$
This immediately implies that $N=M$ and, together with \eqref{e:dord}, that $u_j(x) = v_j(x)$
for $\Leb^n$-a.e.~$x \in \R^n$ for every $j = 1, \ldots , N$.

\subsection{Equality of the total variations}

We know that $\i(u_j)=(\Gamma(u_j),1,\tau_{u_j})$, and the locality
property stated in Lemma~\ref{localtau} allows us to find a Borel
orientation $\tau$ of $\Gamma:=\cup_j\Gamma(u_j)$ with the property
\medskip
\begin{equation}\label{gat2}
\tau=\tau_{u_j}\qquad\text{$\Haus^n$-a.e.~on $\Gamma(u_j)$, for $j=1,\ldots,N$,}
\end{equation}
since by construction the functions $u_j$ satisfy \eqref{e:monot}.
Let us define $\theta(w)$ as the cardinality of the set $\{j\in\{1,\ldots,N\}:\ w\in\Gamma(u_j)\}$;
taking \eqref{gat2} into account, we have then
\medskip
$$
\langle T,\omega\rangle=
\sum_{j=1}^N\langle\i(u_j),\omega\rangle=
\sum_{j=1}^N\int_{\Gamma(u_j)}\langle \tau_{u_j},\omega\rangle\,d\Haus^n
=\int_{\Gamma}\theta\langle \tau,\omega\rangle\,d\Haus^n.
$$
This proves that $T=(\Gamma,\theta,\tau)$. As a consequence
\medskip
$$
\|T\|=\theta\Haus^n\res\Gamma=\sum_{j=1}^N\Haus^n\res\Gamma(u_j)=\sum_{j=1}^N\|\i(u_j)\|.
$$

\medskip
\medskip


\end{document}